\documentclass{casc}
\usepackage[psamsfonts]{amssymb}
\usepackage{latexsym}
\usepackage{enumerate}
\usepackage{amsmath}
\def\span{\mathrm{span}}
\def\chara{\mathrm{char}}
\begin{document}
\title{On the Structure of Cohomology of Hamiltonian $p$-Algebras}
\titlerunning{Cohomology of Hamiltonian $p$-Algebras}
\author{Vladimir V. Kornyak}
\institute{Laboratory of Information Technologies\\
           Joint Institute for Nuclear Research\\
           141980 Dubna, Russia \\
           \email{kornyak@jinr.ru}}
\authorrunning{Vladimir V. Kornyak}
\maketitle
\begin{abstract}
We demonstrate advantages of non-standard grading for computing
cohomology of restricted Hamiltonian and Poisson algebras.
These algebras contain the inner grading element in the properly
defined symmetric grading compatible with the symplectic structure.
Using \emph{modulo $p$} analog of the theorem on the structure of
cohomology of Lie algebra with inner grading element, we show that
all nontrivial cohomology classes are located in the grades which
are the multiples of the characteristic $p$. Besides, this grading
implies another symmetries in the structure of cohomology. These
symmetries are based on the Poincar\'e duality and symmetry with
respect to transpositions of conjugate variables of the symplectic
space. Some results obtained by computer program utilizing these 
peculiarities in the cohomology structure are presented.
\end{abstract}

\section{Introduction}

The ground field is $F$. The Hamiltonian algebra $\mathfrak{h}(n)$ is 
the Lie algebra of vector fields on $F^n$ annihilating the 2-form
\begin{equation}
\omega = \sum_{i=1}^m dx_i\wedge dx_{i+m},\;\text{ where $n = 2m\,$.}
\label{2form}
\end{equation} 
Since this algebra plays a key role in both classical and quantum
physics, it has been widely studied by different methods, in particular,
using the tools of algebraic topology.  Its cohomologies are
important invariants and for the past four decades I.~Gelfand and
his collaborators tried to compute them.  All the computations
performed so far are partial; for an account, see
\cite{Fuks,Fuks85} (the latest result, with Kontsevich,
concerned a version of $\mathfrak{h}(n)$ whose elements were vector
fields with Laurent polynomials as coefficients).

One of the main approaches to deal with the cohomology of a given 
infinite-dimensional algebra is to extract some finite-dimensional
subcompexes based usually on an expertly introduced grading.
Nevertheless, in many important cases, e.g., in investigation of
deformations via cohomology with coefficients in the adjoint module,
this trick does not work.

Another approach is based on construction of
finite-dimensional models of infinite-dimensional algebras.

Any Lie algebra of polynomial vector fields has finite-dimensional
analogs defined over fields of positive characteristic $p$. These
analogs  sometimes possess a special structure; they are called 
\emph{restricted Lie algebras} or \emph{Lie $p$-algebras.} 
These algebras were first systematically studied by Jacobson
in \cite{Jacobson}. The corresponding cohomology theory was first
considered by Hochschild in \cite{Hochschild}. For most general and
mathematically rigorous information on the subject, see
\cite{KostrikinShafarevich6,KostrikinShafarevich9,Strade7,Strade8,%
StradeFarnsteiner,PremetStrade7,PremetStrade9,PremetStrade1}.
Here we consider only few elementary constructions sufficient
for our purposes.

We shall use in what follows the construction called the
\emph{algebra of divided powers}. Let the characteristic of $F$ be
$p>0$ and $\mathbf{x}=\left(x_1,\ldots,x_n\right)$ be a set of indeterminates.
For a multiindex $\mathbf{r} = \left(r_1,\ldots,r_n\right)$, set
\begin{equation}
\mathbf{x}^{(\mathbf{r})} = \prod_{i=1}^{n} x_i^{(r_i)}
= \prod_{i=1}^{n}\frac{x_i^{r_i}}{r_i!}\,.
\label{basis}
\end{equation}
Let $F\left[\mathbf{x}\right]_p =
 \span\left\{\mathbf{x}^{(\mathbf{r})}|r_i < p\right\}$, 
then one can see that $F\left[\mathbf{x}\right]_p$
is a subalgebra of algebra $F\left[\mathbf{x}\right]$ if
$\chara F = p,$ since the multiplication for basis elements of form
(\ref{basis}) is given by the formula
$$
\mathbf{x}^{(\mathbf{r})}\mathbf{x}^{(\mathbf{s})} 
= \binom{\mathbf{r}+\mathbf{s}}{\mathbf{r}}
\mathbf{x}^{(\mathbf{r}+\mathbf{s})}
$$
and $\binom{i+j}{i} = 0 \mod p \ $ for integer $i$ and $j$ such that
$0 \leq i,j < p$ and $i+j \geq p$.

The main result of this paper, presented in Section \ref{mainsection},
is demonstration that the structure of cohomology of Lie $p$-algebra
with the Poisson bracket becomes more clear and easier for computation
if we use the grading compatible with symplectic structure (\ref{2form}).

\section{The Restricted Hamiltonian Algebras}
The elements of $\mathfrak{h}(n)$
--- the
Hamiltonian vector fields --- can be expressed in terms of
\emph{generating functions}.  The elements of the central extension of
the Hamiltonian algebra, the Poisson algebra $\mathfrak{po}(n)$, can
also be naturally described by means of generating function.  

If we express the generating function in terms of divided power
monomials (\ref{basis}), then, in the characteristic $p$, we obtain
the \emph{truncated} Hamiltonian and Poisson algebras denoted in what
follows by $\mathfrak{h}(n)_p$ and $\mathfrak{po}(n)_p$, respectively.
 
Observe that although, for the elements of $F\left[\mathbf{x}\right]_p$, the
degrees of indeterminates are $<p$, the monomials of generating
functions for $\mathfrak{h}(n)_p$ and $\mathfrak{po}(n)_p$ may contain
degrees equal to $p$ since the mapping of the space of generating
functions into the space of vector fields involves differentiations.

We see that
$$
\mathfrak{po}(n)_p = \span\left\{x_1^{(r_1)}\cdots
 x_n^{(r_n)}|0\leq r_i < p\right\}\oplus
\span\left\{x_1^{(p)},\ldots,x_n^{(p)}\right\},
$$
$$
\mathfrak{h}(n)_p = \span\left\{x_1^{(r_1)}\cdots
 x_n^{(r_n)}|0\leq r_i < p, 0 < \sum_{i=0}^{n} r_i
\right\}\oplus
\span\left\{x_1^{(p)},\ldots,x_n^{(p)}\right\},
$$
$$
\dim \mathfrak{po}(n)_p = p^n + n, 
\quad \dim \mathfrak{h}(n)_p = p^n + n - 1.
$$
In what follows, we will mainly restrict our attention to
the Hamiltonian case. The algebra
$\mathfrak{h}(n)_p$ is not semisimple, in other words,
its first cohomology
$H^1\left(\mathfrak{h}(n)_p\right) = 
\left(\mathfrak{h}(n)_p/
\left[\mathfrak{h}(n)_p,
\mathfrak{h}(n)_p\right]\right)^\prime$ is non-trivial.
The ideal 
$$
\mathfrak{h}^{(1)}(n)_p =
\left[\mathfrak{h}(n)_p,
\mathfrak{h}(n)_p\right]
 = \span\left\{x_1^{(r_1)}\cdots x_n^{(r_n)}
 |0\leq r_i < p, 0 < \sum_{i=0}^{n} r_i\right\},
$$
$$ 
\dim \mathfrak{h}^{(1)}(n)_p = p^n - 1
$$
still remains non-semisimple, but the next ideal
$$
\mathfrak{h}^{(2)}(n)_p =
\left[\mathfrak{h}^{(1)}(n)_p,
\mathfrak{h}^{(1)}(n)_p\right]
 = \span\left\{x_1^{(r_1)}\cdots x_n^{(r_n)}
 |0\leq r_i < p, 0 < \sum_{i=0}^{n} r_i < n(p-1)\right\},
$$
$$ 
\dim \mathfrak{h}^{(2)}(n)_p = p^n - 2
$$
is a simple Lie $p$-algebra for $p\neq2$.
The algebra $\mathfrak{h}^{(2)}(n)_p$ is identified
with the Lie $p$-algebra denoted by $H_{m}$  ($m = n/2$) in the well known
list%
\footnote{This list contains four series of algebras: $W_n, S_n, H_n,
K_n$, see \cite{KostrikinShafarevich9}.} of restricted simple Lie
algebras of Cartan type.

\section{A Preliminary Example}
To give some idea of how non-standard grading may be useful, let us
look at an example. We consider here the cohomologies
$H^k_{g}\left(\mathfrak{h}(2)_3\right)$ \footnote{It might seem more
natural to consider as an example some simple ideal instead of
$\mathfrak{h}(2)_3$, but $\mathfrak{h}^{(2)}(2)_3$ is too small to be
illustrative, whereas $\mathfrak{h}^{(2)}(2)_5$ is too large.}
computed by the author's program described in \cite{KornCASC03}. In
Table \ref{standard} the dimensions of
$H^k_{g}\left(\mathfrak{h}(2)_3\right)$ in the \emph{standard} grading
($\deg x_i = 1$ for all $i$) are presented. In this and subsequent
tables an empty box means that $\dim C^k_g=0$, a dot means that 
$\dim C^k_g\neq0$ but $\dim H^k_g=0$.
\begin{table}[h!]
	\caption{$\dim H^k_g\left(\mathfrak{h}(2)_3\right)$ in the standard grading.}
	\label{standard}
	\begin{center}
		\begin{tabular}{l|cccccccccc}
$g\backslash k$&1&2&3&4&5&6&7&8&9&10
\\
\hline
-2&&1&$\cdot$&$\cdot$&1&&&&&
\\				
-1&$\cdot$&$\cdot$&$\cdot$&$\cdot$&$\cdot$&$\cdot$&&&&
\\				
~0&$\cdot$&$\cdot$&4&1&$\cdot$&3&1&&&
\\				
~1&2&$\cdot$&$\cdot$&4&2&$\cdot$&2&2&&
\\				
~2&$\cdot$&$\cdot$&$\cdot$&$\cdot$&$\cdot$&$\cdot$&
$\cdot$&$\cdot$&$\cdot$&
\\				
~3&&2&2&$\cdot$&2&4&$\cdot$&$\cdot$&2&
\\				
~4&&&1&3&$\cdot$&1&4&$\cdot$&$\cdot$&1
\\				
~5&&&&$\cdot$&$\cdot$&$\cdot$&$\cdot$&$\cdot$&$\cdot$&
\\				
~6&&&&&1&$\cdot$&$\cdot$&1&&
\\				
		\end{tabular}
	\end{center}
\end{table}

The picture in Table \ref{standard} does not look instructive at all,
so let us introduce another grading: $\deg x_i = -1$, $\deg x_{i+m} = 1$ for $i\leq m$. We shall call this (of non-Weisfeiler's type) grading \emph{symmetric}. Repeating
computation with symmetric grading we obtain results 
presented in Table
\ref{symmetric}.

\begin{table}[h!]
	\caption{$\dim H^k_g\left(\mathfrak{h}(2)_3\right)$ in the symmetric grading.}
	\label{symmetric}
	\begin{center}
		\begin{tabular}{l|cccccccccc}
$g\backslash k$&1&2&3&4&5&6&7&8&9&10
\\
\hline
-7&&&&$\cdot$&$\cdot$&$\cdot$&&&&
\\				
-6&&&1&1&$\cdot$&1&1&&&
\\				
-5&&$\cdot$&$\cdot$&$\cdot$&$\cdot$&$\cdot$&$\cdot$&
$\cdot$&&
\\				
-4&&$\cdot$&$\cdot$&$\cdot$&$\cdot$&$\cdot$&$\cdot$&
$\cdot$&&
\\				
-3&1&1&1&2&2&2&1&1&1&
\\				
-2&$\cdot$&$\cdot$&$\cdot$&$\cdot$&$\cdot$&$\cdot$&
$\cdot$&$\cdot$&$\cdot$&
\\				
-1&$\cdot$&$\cdot$&$\cdot$&$\cdot$&$\cdot$&$\cdot$&
$\cdot$&$\cdot$&$\cdot$&
\\				
~0&$\cdot$&1&3&2&2&2&3&1&$\cdot$&1
\\				
~1&$\cdot$&$\cdot$&$\cdot$&$\cdot$&$\cdot$&$\cdot$&
$\cdot$&$\cdot$&$\cdot$&
\\				
~2&$\cdot$&$\cdot$&$\cdot$&$\cdot$&$\cdot$&$\cdot$&
$\cdot$&$\cdot$&$\cdot$&
\\				
~3&1&1&1&2&2&2&1&1&1&
\\				
~4&&$\cdot$&$\cdot$&$\cdot$&$\cdot$&$\cdot$&$\cdot$&
$\cdot$&&
\\				
~5&&$\cdot$&$\cdot$&$\cdot$&$\cdot$&$\cdot$&$\cdot$&
$\cdot$&&
\\				
~6&&&1&1&$\cdot$&1&1&&&
\\				
~7&&&&$\cdot$&$\cdot$&$\cdot$&&&&
\\				
		\end{tabular}
	\end{center}
\end{table}

Now the picture is much more attractive. 
The reasons are explained in the next section.
 
\section{The Symmetric Grading}
\label{mainsection}

Let $\mathfrak{a}(n)_p$ denote any
of the algebras $\mathfrak{h}(n)_p$,
 $\mathfrak{h}^{(2)}(n)_p$ or the corresponding Poisson algebras.
Let  $N=\dim \mathfrak{a}(n)_p$.

The standard $\mathbb{Z}$-grading of $\mathfrak{a}(n)_p$ is induced by
prescribing the grades $\deg x_i = 1$ to all $n$ indeterminates. This
grading divides the algebra $\mathfrak{a}(n)_p$ into the homogeneous
degree $i$ subspaces $L_i$:  
\begin{equation}
\mathfrak{a}(n)_p = \left(\ L_{-2}\ \oplus\ \right)\
 L_{-1}\oplus L_0\oplus L_1\oplus \cdots\oplus L_r\,,
 \label{standardgr}
\end{equation}
$\left(\phantom{\frac{1}{1}\!\!}L_{-2} = 0\right.$ for Hamiltonian algebras, 
$r =(p-1)n-3$ for 
$\left.\mathfrak{h}^{(2)}(n)_p\phantom{\frac{1}{1}\!\!}\right)$.

Imposing grading like (\ref{standardgr}) is, perhaps, the only
possibility to extract finite-dimensional subcomplexes
when computing cohomology of infinite-dimensional algebra.

But for finite-dimensional restricted Lie algebra
of vector fields with the Poisson brackets we can use the
symmetric grading $\deg x_i = -1,\ \deg x_{i+m} = 1$ reflecting 
the fact that $x_i$ and $x_{i+m}$ are 
conjugate indeterminates for symplectic structure (\ref{2form}). 
This grading induces the following decomposition:
\begin{equation}
\mathfrak{a}(n)_p = L_{-r}\oplus \cdots\oplus L_{-1}\oplus L_0\oplus
L_1\oplus \cdots\oplus L_r\,,
 \label{symmetricgr}
\end{equation}
$\left(r =\frac{(p-1)n}{2}\ \right.$ for 
$\left.\mathfrak{h}^{(2)}(n)_p
\phantom{\frac{(p-1)n}{2}}\hspace{-25pt}\right)$.
The advantages of grading (\ref{symmetricgr}) 
are clear from the following propositions 1--3.

\textbf{Proposition 1.}
\emph{Cohomologies in subcomplexes with opposite grades are isomorphic:}
\begin{equation}
H^k_{g}\left(\mathfrak{a}(n)_p\right)\cong
H^k_{-g}\left(\mathfrak{a}(n)_p\right)\,.
\label{prop1}
\end{equation}
\emph{Proof:} The isomorphism follows from the symmetry of all
 constructions with respect to
transpositions of the conjugate indeterminates $x_i$ and $x_{i+m}$.
\hfill$\square$\\[2pt]
\emph{Comment:} Turning to Table \ref{symmetric},
we see that with Proposition 1 it suffices to compute only
one half (upper or lower) of the table.

\textbf{Proposition 2.}
\emph{Cohomology of degree $k$ in a given grade $g$ is isomorphic
 to the homology of degree $N-k$ in the same grade $g$:}
\begin{equation}
H^k_{g}\left(\mathfrak{a}(n)_p\right)\cong
\left(H^{N-k}_{g}\left(\mathfrak{a}(n)_p\right)\right)^{\prime}
=H_{N-k,g}(\mathfrak{a}(n)_p)\,.
\label{prop2}
\end{equation}
\emph{Proof:} The proposition follows immediately from 
the \emph{Poincar\'e duality.} The algebra $\mathfrak{a}(n)_p$ 
is \emph{unitary}  \cite{Fuks}, i.e., $H_N\left(\mathfrak{a}(n)_p\right)\neq 0.$ 
Since in the symmetric grading $C_N\left(\mathfrak{a}(n)_p\right)$ has zero grade,
the Poincar\'e dual to the 
$H^k_{g}\left(\mathfrak{a}(n)_p\right)$ is 
$H_{N-k,-g}(\mathfrak{a}(n)_p)$, and using (\ref{prop1}) to change 
sign of grade we obtain (\ref{prop2}).
\hfill$\square$\\[2pt]
\emph{Comment:} Proposition 2 also reduces the computation to one half
 (left or right) of the table. Note also, that we, 
 slightly violating the symmetry of the picture, did not include in Table 
\ref{symmetric} one-dimensional $H^0_{0}\left(\mathfrak{h}(2)_3\right)$ 
which is dual to the $H_{10,0}\left(\mathfrak{h}(2)_3\right)$ and is 
isomorphic to the ground field $F$.

\textbf{Proposition 3.} \emph{For the symmetric grading all non-trivial cohomologies
$H^k_{g}\left(\mathfrak{a}(n)_p\right)$ lie in the grades $g = \pm pj$ for $j$ integer.}\\[5pt]
\emph{Proof:} This proposition is, in fact, ``modulo $p$'' version 
of standard theorem \cite{Fuks} on the \emph{inner grading element}. 
One should only reformulate slightly the theorem and repeat the proof 
replacing the arithmetic in zero characteristic by the modular arithmetic.
The inner grading element for $\mathfrak{a}(n)_p$ in symmetric grading
is represented by the generating function $\sum^{m}_{i=1}x_ix_{i+m}.$
\hfill$\square$\\[2pt]
\emph{Comment:} Proposition 3 manifests in Table \ref{symmetric} in 
the fact that non-trivial cohomologies are located in the rows $-6$, $-3$, $0$,
$3$, $6$ only.

From the computational point of view, the combination of Propositions
1--3 reduces the computation of the cohomology roughly by the factor 
$1/4p.$ Turning to Table \ref{symmetric} we see that there is 108 
$(g,k)$-subcomlexes with non-empty cochain spaces, and using 
the above propositions it suffices to process only
13 of them. Note also, that in the standard grading
we have (see Table \ref{standard}) only 60 non-empty
$(g,k)$-subcomlexes, i. e., the average size of subcomplex 
in the symmetric grading is smaller than in the standard grading (and,
hence, on the average, the subcomplexes should be easier for computation).

\section{Computation of $H^*_*\left(\mathfrak{h}^{(2)}(2)_5\right)$}

Since $\dim \mathfrak{h}^{(2)}(2)_5=23$, 
the total (i.e., for all degrees and grades) dimension of the
cochain space is
$\dim C^*_*\left(\mathfrak{h}^{(2)}(2)_5\right) = 8388608.$ 
Thus, the full computation of the cohomology
is a rather difficult task. Nevertheless we have completed the task. Using the above listed symmetries 
in the structure of cohomology it suffices to carry out computations 
for degrees $1\leq k\leq 11$ 
in the grades 0, 5, 10, 15, 20.
The data presented in Table \ref{h5} were obtained on a 1133MHz Pentium III PC with 512Mb under the Windows XP.
The computation of $H^{11}_0\left(\mathfrak{h}^{(2)}(2)_5\right)$ took 22 h 29 min.
It is the most difficult subtask. All other $(g,k)$-boxes of the table take much smaller time of calculations. 

\begin{table}[h!]
	\caption{$\dim H^k_g\left(\mathfrak{h}^{(2)}(2)_5\right)$.}
	\label{h5}
	\begin{center}
		\begin{tabular}{l|ccccccccccc}
$g~\backslash~ k$
&~1~&~2~&~3~&~4~&~5~&~6~&~7~&~8~&~9~&~10~&~11~
\\[2pt]
\hline\\[-5pt]
$~~~0$&$\cdot$&1&1&4&4&8&12&9&18&14&30
\\[2pt]				
$\pm~5~$&&1&1&3&2&6&9&8&15&14&25
\\[2pt]				
$\pm10~$&&&$\cdot$&3&1&3&6&4&9&7&17
\\[2pt]				
$\pm15~$&&&&&&2&1&1&4&3&7
\\[2pt]				
$\pm20~$&&&&&&&&&&1&3
\\				
		\end{tabular}
	\end{center}
\end{table}

\section{Conclusions}
As is clear, the full computations for $p > 5$ are hardly possible on the present-day computers,
e.g., $\dim \mathfrak{h}^{(2)}(2)_7 = 47$ and
$\dim C^*_*\left(\mathfrak{h}^{(2)}(2)_7\right) = 140737488355328$.
Applying different tricks and efforts we can, probably, reduce the
last value by several orders, but this is quite insufficient.

On the other hand, we can try to derive some hints 
about the structure of cohomology of classical Hamiltonian algebra 
analyzing the small degree cohomologies of $p$-analogs of the classical algebra.
We have performed some computations for $p = 7, 11, 13$
and it seems that the classical cohomology classes
are among the zero grade cocycles of $p$-algebras.
Of course, for any finite $p$,~ zero grade series of
cocycles contain multiplicative combinations
of cocycles in grades $\pm p, \pm 2p,\ldots,$ etc.
and these should be separated from the classical cocycles.

It seems also that the sequences of dimensions
in the grades $jp$ with growing $p$ tend
to some stable sequences depending only on the number $j$
(not on $p$) but these observations are too preliminary to be 
discussed seriously at present.

\section*{Acknowledgments}
I would like to thank D. Leites for raising the problem.
This work was supported in part by the 
grants 04-01-00784 from the Russian Foundation for Basic Research and
2339.2003.2 from the Russian Ministry of Industry, Science and
Technologies.

\end{document}